\documentclass[12pt,a4paper]{article}%
\usepackage{amsmath}
\usepackage{amsfonts}
\usepackage{amssymb}
\usepackage{graphicx}%
\numberwithin{equation}{section}
\newtheorem{theorem}{Theorem}

\newtheorem{corollary}{Corollary}

\newtheorem{definition}[theorem]{Definition}

\newtheorem{proposition}{Proposition}
\newtheorem{remark}{Remark}

\newenvironment{proof}[1][Proof]{\noindent\textbf{#1.} }{\ \rule{0.5em}{0.5em}}

\def \f {\frac}

\def \O {\Omega}
\def \R {\mathbb R}

\begin{document}

\title{Steady free convection in a bounded and saturated porous medium}
\author{Samir AKESBI$\dag$, Bernard BRIGHI$\ddagger$  and Jean-David HOERNEL$\sharp$ }
\date{}
\maketitle

\begin{center}
Universit\'e de Haute-Alsace, Laboratoire de Math\'ematiques, Informatique et Applications
\vskip 0,1cm
4 rue des fr\`eres Lumi\`ere, 68093 MULHOUSE (France)
\vskip 1,5cm
\end{center}

\begin{abstract} In this paper we are interested with a strongly coupled system of partial differential equations that modelizes free convection in a two-dimensional bounded domain filled with a fluid saturated porous medium. This model is inspired by the one of free convection near a semi-infinite impermeable vertical flat plate embedded in a fluid saturated porous medium. We establish the existence and uniqueness  of the solution for small data in some unusual spaces.
\end{abstract}

\footnotetext{AMS 2000 Subject Classification: 34B15, 34B40, 35Q35, 76R10, 76S05.}
\footnotetext{Key words and phrases: Free convection, porous medium, coupled pdes.}
\footnotetext{$\dag$ s.akesbi@uha.fr $\ddagger$ b.brighi@uha.fr $\sharp$ j-d.hoernel@wanadoo.fr}

\section{Introduction}
In the literature, many papers about free convection in fluid saturated porous media study the case of the semi-infinite vertical flat plate in the framework of boundary layer approximations. This approach allows to introduce similarity variables to reduce the whole system of partial differential equations into one single ordinary differential equation of the third order with appropriate boundary values. This two points boundary value problem can be studied using a shooting method or an auxiliary dynamical system either in the case of prescribed temperature or in the case of prescribed heat flux along the plate. 

In this article we first present the derivation of the equations, show how the boundary layer approximation leads to the two points boundary value problem and the similarity solutions, then we rewrite the full problem of free convection in a two-dimensional bounded domain filled with a fluid saturated porous medium. This new model, written in terms of stream function and temperature, consists in two strongly coupled partial differential equations. We establish the existence and uniqueness of its solution for small data.
\bigskip

\section{The semi-infinite vertical flat plate case}
Let us consider a semi-infinite vertical permeable or impermeable flat plate embedded in a fluid saturated porous medium at the ambient temperature $T_{\infty}$, and a rectangular Cartesian co-ordinates system with the origin fixed at the leading edge of the vertical plate, the $x$-axis directed upward along the plate and the $y$-axis normal to it. If we suppose that the porous medium is homogeneous and isotropic, that all the properties of the fluid and the porous medium are constants and that the fluid is incompressible and follows the Darcy-Boussinesq law we obtain the following governing equations 
$$\frac{\partial u}{\partial x}+\frac{\partial v}{\partial y}=0,$$
$$u=-\frac{k}{\mu}\left(\frac{\partial p}{\partial x}+\rho g \right),$$
$$v=-\frac{k}{\mu}\frac{\partial p}{\partial y},$$
$$u\frac{\partial T}{\partial x}+v\frac{\partial T}{\partial y}=
\lambda \left(\frac{\partial^2 T}{\partial x^2}
+\frac{\partial^2 T}{\partial y^2}\right ),$$
$$\rho=\rho_{\infty}(1-\beta(T-T_{\infty}))$$
in which $u$ and $v$ are the Darcy velocities in the $x$ and $y$ directions, $\rho$, $\mu$ and
$\beta$ are the density, viscosity and thermal expansion coefficient of the fluid, $k$ is the
permeability of the saturated porous medium, $\lambda$ is its thermal diffusivity, $p$
is the pressure, $T$ the temperature and $g$ the acceleration of the gravity. The subscript 
$\infty$ is used for values taken far from the plate. In our system of co-ordinates there are two main interesting sets of boundary conditions along the plate. 

First, the temperature is prescribed on the wall that gives
\begin{equation}
v(x,0)=\omega x^\frac{m-1}{2}, \quad T(x,0)=T_w(x)=T_\infty+Ax^m \label{b1}
\end{equation}
with $m\in\mathbb{R}$ and $A>0$, see \cite{pop1}, \cite{cheng}, \cite{ene}, \cite{ing} and \cite{merk}.

Secondly, the heat flux is prescribed along the plate that leads to
\begin{equation}
v(x,0)=\omega x^\frac{m-1}{3},
\quad \frac{\partial T}{\partial y}(x,0)=-x^m \label{b2}
\end{equation}
with $m\in\mathbb{R}$, see \cite{heat_flux} and \cite{pop}.

The parameter $\omega \in \mathbb{R}$ is the mass transfer coefficient. For an impermeable wall we have $\omega=0$, and for a permeable wall, $\omega<0$ corresponds to  fluid suction and $\omega>0$ to fluid injection. The boundary conditions far from the plate are the same in both cases (\ref{b1}) and (\ref{b2})
$$u(x,\infty)=0,\quad T(x,\infty)=T_{\infty}.$$
If we introduce the stream function $\Psi$ such that
$$u=\frac{\partial \Psi}{\partial y}, \quad v=-\frac{\partial \Psi}{\partial x}$$
we obtain the system in which it remains only $\Psi$ and $T$
\begin{equation}
\frac{\partial^2 \Psi}{\partial x^2}+\frac{\partial^2 \Psi}{\partial y^2}
=\frac{\rho_\infty \beta g k}{\mu}\frac{\partial T}{\partial y}, \label{p11}
\end{equation}
\begin{equation}
\lambda \left(\frac{\partial^2 T}{\partial x^2}+
\frac{\partial^2 T}{\partial y^2} \right)=\frac{\partial T}{\partial x}
\frac{\partial \Psi}{\partial y}-\frac{\partial T}{\partial y}
\frac{\partial \Psi}{\partial x}. \label{p12}
\end{equation}
Along the wall, the boundary conditions (\ref{b1}) become
\begin{equation}
\frac{\partial \Psi}{\partial x}(x,0)=-\omega x^\frac{m-1}{2}, \quad T(x,0)=T_w(x)=T_\infty+Ax^m \label{bp1}
\end{equation}
and (\ref{b2}) becomes
\begin{equation}
\frac{\partial \Psi}{\partial x}(x,0)=-\omega x^\frac{m-1}{3},
\quad \frac{\partial T}{\partial y}(x,0)=-x^m.\label{bp2}
\end{equation}
The boundary conditions far from the plate become
\begin{equation}
\frac{\partial \Psi}{\partial y}(x,\infty)=0,\quad T(x,\infty)=T_{\infty}.\label{bf}
\end{equation}
We will start from the equations (\ref{p11})-(\ref{p12}) subjected to the boundary conditions (\ref{bp1}) and (\ref{bf}) with $\omega=0$ to write a new model, settled in a two-dimensional bounded domain, that we will study in the rest of this paper.

Before doing this, let us say a few words about the similarity solutions. Assuming that convection takes place in a thin layer around the plate, we obtain the boundary layer approximation
\begin{equation}
\frac{\partial^2 \Psi}{\partial y^2}
=\frac{\rho_\infty \beta g k}{\mu}\frac{\partial T}{\partial y}, \label{epsi}
\end{equation}
\begin{equation}
\frac{\partial^2 T}{\partial y^2}=\frac{1}{\lambda}\left (
\frac{\partial T}{\partial x}
\frac{\partial \Psi}{\partial y}-\frac{\partial T}{\partial y}
\frac{\partial \Psi}{\partial x} \right ) \label{T} 
\end{equation}
with the same boundary conditions (\ref{bp1}) or (\ref{bp2}) and (\ref{bf}) as before. 

For the case of prescribed heat, introducing the new dimensionless similarity variables
$$t=(Ra_{x})^\frac{1}{2}\frac{y}{x},\quad
\Psi(x,y)=\lambda (Ra_{x})^\frac{1}{2}f(t),\quad
T(x,y)=(T_w(x)-T_\infty)\theta(t)+T_\infty$$
with 
$$Ra_x=\frac{\rho_\infty \beta g k(T_w(x)-T_\infty)x}{\mu\lambda}$$
the local Rayleigh number, equations (\ref{epsi}) and (\ref{T}) with the boundary conditions (\ref{bp1}) and (\ref{bf}) leads to the third order ordinary differential equations
$$f'''+\frac{m+1}{2}ff''-mf'^2=0$$
on $[0,\infty)$ subjected to
$$f(0)=-\gamma, \quad f'(0)=1\quad \text{and}\quad f'(\infty)=0$$
where
$$\gamma=\frac{2\omega}{m+1}\sqrt{\frac{\mu}{\rho_\infty \beta g k A \lambda}}.$$
One can find explicit solutions of this problem for some particular values of $\gamma$ or $m$ in \cite{brighicr}, \cite{brighi02}, \cite{brighi01}, \cite{crane}, \cite{gup},  \cite{ing}, \cite{mag} and \cite{stu}. For mathematical results about existence, nonexistence, uniqueness, nonuniqueness and asymptotic behavior, see \cite{banks1}, \cite{brighicr}, \cite{brighi02} and \cite{ing} for $\gamma=0$, and  \cite{brighi01},  \cite{equiv}, \cite{brighisari}, \cite{guedda} and  \cite{guedda1} for the general case. Numerical  investigations can be found in \cite{banks1}, \cite{brighi04}, \cite{pop1}, \cite{cheng}, \cite{ing}, \cite{mag} and \cite{wood}.

In the case of prescribed heat flux, we introduce the new dimensionless similarity variables
$$t=3^{-\frac{1}{3}}R_{a}^\frac{1}{3}x^\frac{m-1}{3}y,\quad
\Psi(x,y)=3^{\frac{2}{3}} R_{a}^\frac{1}{3}\lambda x^\frac{m+2}{3}f(t),\quad
T(x,y)=3^{\frac{1}{3}}R_{a}^{-\frac{1}{3}}x^\frac{2m+1}{3}\theta(t)+T_{\infty}$$
and the Rayleigh number
$$R_{a}=\frac{\rho_{\infty}\beta g k}{\mu\lambda}.$$
Then, equations (\ref{epsi}) and (\ref{T}) with the boundary conditions (\ref{bp2})-(\ref{bf}) give
$$f'''+(m+2)ff''-(2m+1)f'^2=0$$
and
$$f(0)=-\gamma, \quad f''(0)=-1\quad \text{and}\quad f'(\infty)=0$$
where
$$\gamma=\frac{3^\frac{1}{3}R_a^{-\frac{1}{3}}\omega}{\lambda (m+2)}.$$
The study of existence, uniqueness and qualitative properties of the solutions of this problem is made in \cite{heat_flux}. For a survey of the two cases, see \cite{gaeta}. This equation
can also be found in industrial processes such as boundary layer flow adjacent to stretching walls (see \cite{banks1}, \cite{banks}, \cite{crane}, \cite{gup}, \cite{mag}) or excitation of liquid metals in a high-frequency magnetic field (see \cite{mof}).

One particular case of the two previous equations is the Blasius equation $f'''+ff''=0$ introduced in \cite{bla} and studied, for example, in \cite{brighi03},  \cite{coppel} and \cite{hart}.

The case of mixed convection $f'''+ff''+mf'(1-f')=0$ with $m\in\mathbb R$ is interesting too and results about it can be found in \cite{aly}, \cite{aml}, \cite{guedda2} and \cite{nazar}. The Falkner-Skan equation $f'''+ff''+m(1-f'^2)=0$ with $m\in\mathbb R$ is in the same family of problems, see \cite{coppel}, \cite{falk}, \cite{hart}, \cite{ish1}, \cite{wang}, \cite{yang} and \cite{yang2} for results about it.

New results about the more general equation $f'''+ff''+g(f')=0$ for some given function $g$ can be found in \cite{jde}, see also \cite{utz}.

\section{A model problem in a bounded domain}
Let $\Omega\subset \R^2$ be a simply connected, bounded lipschitz domain whose boundary 
$\Gamma=\partial \O$ is divided in two connected parts $\Gamma_1$ and $\Gamma_2$ such that 
$$\overline{\Gamma}_{1} \cup \overline{\Gamma}_{2}=\Gamma \text{ and } \Gamma_{1} \cap \Gamma_{2}=\emptyset.$$
\begin{center}
\includegraphics{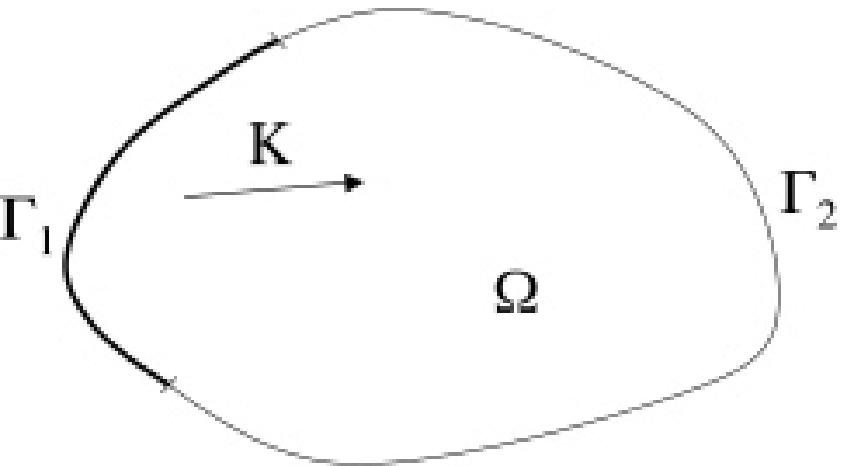}
\end{center}
We start from the previous equations (\ref{p11})-(\ref{p12}) in terms of the stream function $\Psi$ and the temperature $T$ with $K=\left (0,\frac{\rho_\infty \beta g k}{\mu}\right )$, and assuming that $\Gamma_1$ is impermeable and that the temperature $T_w\geq  0$ is known on the whole boundary $\Gamma$, we modify the equation (\ref{p11}) by setting $K(x)=(k_1(x),k_2(x)) \in \R^{2}$ with $0<\|K\|_\infty <\infty$. Then, we obtain the following new problem in the bounded domain $\O$, which consists in finding $(\Psi,T)$
$$\Psi:\O \rightarrow \R$$ 
$$T:\O \rightarrow \R$$
verifying the equations in $\O$
\begin{align}
\Delta \Psi&=K.\nabla T, \label{eq1}\\
\lambda \Delta T&= \nabla T. (\nabla \Psi)^\bot, \label{eq2}
\end{align}
the boundary conditions on $\Gamma$ for $\Psi$
\begin{equation}
\Psi=0 \text{ on } \Gamma_{1} \quad \text{and} \quad \frac{\partial \Psi}{\partial n}=0 \text{ on } \Gamma_{2} \label{b1}
\end{equation}
and the boundary conditions on $\Gamma$ for $T$
\begin{equation}
T=T_{w} \text{ on } \Gamma \label{b2}
\end{equation}
where $\lambda \in \R^{+*}$ and for all $x=(u,v)\in \O$, let $x^\bot=(v,-u)$.

\bigskip

\subsection{Preliminary results}

Let us assume that $T_w \in H^\f{1}{2}(\Gamma)$ and let $\Theta$ be the unique function in $H^1(\O)$ verifying 
\begin{align}
\Delta \Theta&=0 \quad \text {in $\O$,} \label{pt1} \\ 
\Theta&=T_{w}  \quad \text{on  $\Gamma$.} \label{pt2} 
\end{align}
In the following we will need that $\nabla \Theta \in L^\infty(\O)$, thus we will suppose that it holds (it is the case if $T_w\in H^\f{5}{2}(\Gamma)$ for example). 

If $(\Psi,T)$ is a solution of $(\ref{eq1})$-$(\ref{b2})$ and if we set $H=T-\Theta$, then $(\Psi,H)$ is a solution of
\begin{align}
\Delta \Psi&=K.\nabla H + K.\nabla \Theta, \label{new1}\\
\lambda \Delta H &= \nabla H. (\nabla \Psi)^\bot + \nabla \Theta. (\nabla \Psi)^\bot
\label{new2}
\end{align}
in the domain $\Omega$ with the boundary conditions for $\Psi$
\begin{equation}
\Psi=0 \text{ on } \Gamma_{1} \quad \text{and} \quad \frac{\partial \Psi}{\partial n}=0 \text{ on } 
\Gamma_{2} \label{bh1}
\end{equation}
and the boundary conditions for $H$
\begin{equation}
H=0 \text{ on } \Gamma. \label{bh2}
\end{equation}
Conversly, it is clear that if $(\Psi,H)$ is a solution of $(\ref{new1})$-$(\ref{bh2})$ then $(\Psi,T):=(\Psi,H+\Theta)$ is a solution of $(\ref{eq1})$-$(\ref{b2})$.

In the following we set $\|.\|_{L^1(\O)}=\|.\|_{1}$, $\|.\|_{L^2(\O)}=\|.\|_{2}$, $\|.\|_{L^\infty(\O)}=\|.\|_{\infty}$ and
$$(u,v)=\int_\O uvdx.$$

\begin{definition}
For $u\in L^\infty(\O)$, $v\in H^1_0(\O)$ and $w \in H^1(\O)$ let 
$$a(u,v,w)=(u\nabla v,(\nabla w)^\bot)_{L^2(\O),L^2(\O)}.$$
\end{definition}
\begin{remark}The trilinear form $a$ is well defined because for $u\in L^\infty(\O)$, $v\in H^1_0(\O)$ 
and $w \in H^1(\O)$ we have
$$|a(u,v,w)|\leq \|u\|_\infty \|\nabla v\|_2\|\nabla w\|_2.$$
\end{remark}

\begin{proposition} \label{a} For $u\in H^1_0(\O)\cap L^\infty(\O)$ and $v \in H^1(\O)$ we have
\begin{equation}
a(u,u,v)=0 \label{int0}.
\end{equation}
\end{proposition}
\begin{proof}
First, let us notice that if $u\in H^1_0(\O)\cap L^\infty(\O)$ then $u^2 \in H^1_0(\O)$ and $\nabla(u^2)=2u\nabla u$. Hence
\begin{align*}
a(u,u,v)&=(u\nabla u,(\nabla v)^\bot)_{L^2(\O),L^2(\O)}\\
&=\f{1}{2}(\nabla u^2,(\nabla v)^\bot)_{L^2(\O),L^2(\O)}\\
&=-\f{1}{2}(\text{div}((\nabla v)^\bot),u^2)_{H^{-1}(\O),H^1_0(\O)}\\
&=0
\end{align*}
because $u=0$ on $\Gamma$ and $\text{div}((\nabla v)^\bot)=0$ in $H^{-1}(\O)$.
\end{proof}

\begin{remark}\label{-a} For $u, v\in H^1_0(\O)\cap L^\infty(\O)$ and $w \in H^1(\O)$ we have
\begin{equation}
a(u,v,w)=-a(v,u,w). \label{int1}
\end{equation}
\end{remark}

\subsection{A priori estimates}

Let 
\begin{align*}
W_\Psi=\left \{u\ | \ u\in H^1(\O) \text{ and } u=0 \text{ on } \Gamma_1 \right \}
\end{align*}
and
$$W_H=H^1_0(\O)\cap L^\infty(\O).$$
The spaces $W_\Psi$ and $W_H$ are equipped with the norms $\|.\|_{W_\Psi}$ and  $\|.\|_{W_H}$ defined by
$$\|u\|_{W_\Psi}=\|\nabla u\|_2 \quad \text{and} \quad \|u\|^2_{W_H}=\|u\|^2_{\infty}+\|\nabla u\|^2_{2}.$$
In the following we will use the notation $C$ for the Poincar\'e's constant of $\O$.
\begin{definition}
We will call $(\Psi,H)\in W_\Psi \times W_H$ a weak solution of the problem $(\ref{new1})$-$(\ref{bh2})$ if and only if we have
\begin{align}
(\nabla \Psi,\nabla u)+(K.\nabla H,u)+(K.\nabla \Theta,u)&=0, \label{v1}\\
\lambda (\nabla H,\nabla v)+a(v,H,\Psi)+a(v,\Theta,\Psi)&=0\label{v2}
\end{align}
for all $u \in W_\Psi$ and $v \in W_H$.
\end{definition}

\begin{proposition}
Let $(\Psi,H)\in W_\Psi \times W_H$ be a solution of the problem $(\ref{v1})$-$(\ref{v2})$ and $T=H+\Theta$, then
\begin{equation}
\inf_{\Gamma} T_w \leq T \leq \sup_{\Gamma} T_w.
\end{equation}
\end{proposition}
\begin{proof}
Set $l=\sup_{\Gamma} T_w $ and $T^+=\sup(T-l,0)$. As $T^+\in W_H$, using $(\ref{v2})$ with $v=T^+$ and noticing that $(\nabla \Theta,\nabla T^+)=0$ because $\Delta \Theta=0$, leads to
$$\lambda (\nabla T,\nabla T^+)+a(T^+,T,\Psi)=0.$$
Using the facts that $\lambda (\nabla T,\nabla T^+)=\lambda (\nabla T^+,\nabla T^+)$ and 
$a(T^+,T,\Psi)=a(T^+,T^+,\Psi)=0$ by proposition $\ref{a}$ we obtain that
$$\| \nabla T^+\|_2=0$$
and as $T^+ \in H^1_0(\O)$ we have $T^+=0$ on $\O$. We proceed in the same way with 
$l'=\inf_{\Gamma} T_w $ and $T^-=\inf(T-l',0)$ for the other inequality.
\end{proof}

\bigskip
\begin{proposition}\label{borne}
Let $(\Psi,H)\in W_\Psi \times W_H$ be a solution of the problem $(\ref{v1})$-$(\ref{v2})$, then for $\|\nabla \Theta\|_\infty <\f{\lambda}{2C^2\|K\|_\infty}$ we have
$$\|\nabla \Psi\|_2 \leq 2C\|K\|_\infty \|\nabla \Theta\|_2\quad \text{and}\quad \|\nabla H\|_2 \leq \|\nabla \Theta\|_2.$$
\end{proposition}
\begin{proof}
Taking $u=\Psi$ in $(\ref{v1})$ and using Poincar\'e's inequality we obtain
\begin{align*}
\|\nabla \Psi\|_2^2 \leq& |(K.\nabla H,\Psi)|+|(K.\nabla \Theta,\Psi)| \\
\leq& \|K\|_\infty \left ( \|\nabla H\|_2+\|\nabla \Theta\|_2 \right ) \|\Psi\|_2 \\
\leq& C\|K\|_\infty \left ( \|\nabla H\|_2+\|\nabla \Theta\|_2 \right ) \|\nabla \Psi\|_2
\end{align*}
and
\begin{equation}
\|\nabla \Psi\|_2 \leq C\|K\|_\infty\left ( \|\nabla H\|_2+\|\nabla \Theta\|_2 \right ). \label{psi}
\end{equation}
Taking $v=H$ in $(\ref{v2})$ leads to
$$\lambda (\nabla H,\nabla H)+a(H,H,\Psi)+a(H,\Theta,\Psi)=0.$$
Then, by proposition $\ref{a}$ we have
\begin{align*}
\lambda \|\nabla H\|_2^2 &\leq |a(H,\Theta,\Psi)| \\
&\leq \|\nabla \Theta\|_\infty \|H\|_2 \|\nabla \Psi\|_2 \\
&\leq C\|\nabla \Theta\|_\infty \|\nabla H\|_2 \|\nabla \Psi\|_2
\end{align*}
using Poincar\'e's inequality and
\begin{equation}
\|\nabla H\|_2 \leq \f{C}{\lambda}\|\nabla \Theta\|_\infty \|\nabla \Psi\|_2. \label{h}
\end{equation}
Then, combining $(\ref{psi})$ and $(\ref{h})$ leads to
\begin{equation*}
\|\nabla \Psi\|_2 \leq C\|K\|_\infty\|\nabla \Theta\|_2 +\f{C^2\|K\|_\infty}{\lambda}\|\nabla \Theta\|_\infty \|\nabla \Psi\|_2.
\end{equation*}
Thus 
\begin{equation*}
\left (1-\f{C^2\|K\|_\infty}{\lambda}\|\nabla \Theta\|_\infty\right )\|\nabla \Psi\|_2 \leq C\|K\|_\infty\|\nabla \Theta\|_2
\end{equation*}
and as $\f{C^2\|K\|_\infty}{\lambda}\|\nabla \Theta\|_\infty <1/2$ we have
\begin{equation*}
\|\nabla \Psi\|_2 \leq 2C\|K\|_\infty\|\nabla \Theta\|_2.
\end{equation*}
Using this new inequality in $(\ref{h})$, we obtain
\begin{equation*}
\|\nabla H\|_2 \leq \|\nabla \Theta\|_2.
\end{equation*}
\end{proof}
\begin{remark}
As $$\|\nabla \Theta\|_2 \leq (\text{\rm mes } \O)^\f{1}{2}\|\nabla \Theta\|_\infty
\quad \text{and}\quad \|\nabla \Theta\|_\infty < \f{\lambda}{2C^2\|K\|_\infty}$$
we can rewrite the previous result as
$$\|\nabla \Psi\|_2 \leq \f{ \lambda }{C}(\text{\rm mes }\O)^\f{1}{2}
\quad \text{and}\quad \|\nabla H\|_2 \leq \f{\lambda }{2C^2\|K\|_\infty}(\text{\rm mes }\O)^\f{1}{2}.$$
\end{remark}

\subsection{Main  results}

\begin{theorem} \label{unicite}
Let $M=\sup_\Gamma{T_w}$. If $MC\|K\|_\infty<\lambda$, then the problem $(\ref{new1})$-$(\ref{bh2})$ admits at most one weak solution $(\Psi,H)$ in $W_\Psi \times W_H$.
\end{theorem}
\begin{proof}
Let $(\Psi_1,H_1)$ and $(\Psi_2,H_2)$ be two solutions of $(\ref{new1})$-$(\ref{bh2})$. Setting $\bar H=H_1-H_2$ and $\bar \Psi=\Psi_1-\Psi_2$ we obtain
\begin{align*}
(\nabla \bar \Psi,\nabla u)+(K.\nabla \bar H,u)&=0, \\
\lambda (\nabla \bar H,\nabla v)+a(v,H_1,\Psi_1)-a(v,H_2,\Psi_2)+a(v,\Theta,\bar \Psi)&=0\ 
\end{align*}
for $u \in W_\Psi$ and $v \in W_H$. Choosing $u=\bar \Psi$ and $v=\bar H$ leads to
\begin{align}
(\nabla \bar \Psi,\nabla \bar \Psi)+(K.\nabla \bar H,\bar \Psi)&=0, \label{u1}\\
\lambda (\nabla \bar H,\nabla \bar H)+a(\bar H,H_1,\Psi_1)-a(\bar H,H_2,\Psi_2)
+a(\bar H,\Theta,\bar \Psi)&=0. \label{u2}
\end{align}
From equation $(\ref{u1})$ we deduce that
\begin{equation}
\|\nabla \bar \Psi\|_2\leq C\|K\|_\infty \|\nabla \bar H\|_2. \label{u3}
\end{equation}
Let us compute
\begin{align*}
a(\bar H,H_1,\Psi_1)-a(\bar H,H_2,\Psi_2) &=-a(H_2,H_1,\Psi_1)-a(H_1,H_2,\Psi_2)\\
&=a(H_1,H_2,\Psi_1)-a(H_1,H_2,\Psi_2)\\
&=a(H_1,H_2,\bar \Psi)\\
&=a(\bar H,H_1,\bar \Psi).
\end{align*}
Thus, using now equation $(\ref{u2})$ we get
$$\lambda (\nabla \bar H,\nabla \bar H)+a(\bar H,H_1+\Theta,\bar \Psi)=0$$
and
\begin{align*}
\lambda \|\nabla \bar H\|_2^2 &\leq |a(\bar H,H_1+\Theta,\bar \Psi)| \\
&\leq |a(T_1,\bar H,\bar \Psi)| \\
&\leq \|T_1\|_\infty \|\nabla \bar H\|_2 \|\nabla \bar \Psi\|_2 \\
&\leq M \|\nabla \bar H\|_2 \|\nabla \bar \Psi\|_2
\end{align*}
with $M=\sup_\Gamma{T_w}$. Therefore
\begin{eqnarray*}
\|\nabla \bar H\|_2&\leq & \frac{M}{\lambda}\|\nabla \bar \Psi\|_2
\end{eqnarray*}
and using $(\ref{u3})$ we have
\begin{eqnarray*}
\|\nabla \bar H\|_2&\leq & \frac{MC\|K\|_\infty }{\lambda}\|\nabla \bar H\|_2.
\end{eqnarray*}
Choosing $\frac{MC\|K\|_\infty }{\lambda}<1$ we obtain $ \|\nabla \bar H\|_2=0$ and $ \|\nabla \bar \Psi \|_2=0$. This complete the proof.
\end{proof}


\bigskip
In the following Theorem, we prove the existence of a strong solution $(\Psi,H)$ of the problem $(\ref{new1})$-$(\ref{bh2})$ under some hypothesis on the data. To this aim, let us define the spaces
 $$\tilde W_\Psi=\left \{u\ | \ u\in H^2(\O), u=0 \text{ on } \Gamma_1 \text{ and } \frac{\partial u}{\partial n}=0 \text{ on } \Gamma_2 \right \}$$
 and $$\tilde W_H=H^1_0(\O)\cap H^2(\Omega).$$
 These spaces are equipped with the following norms
\begin{align*}
\|u\|^2_{\tilde W_\Psi}&=\|\nabla u\|^2_{H^1(\O)},\\
\|v\|^2_{\tilde W_H}&=\|\nabla v\|^2_{H^1(\O)},\\
\| (u,v)\|^2_{\tilde W_\Psi \times \tilde W_H}&=\|u\|^2_{\tilde W_\Psi}+\|v\|^2_{\tilde W_H}
\end{align*}
and
$$\|(u,v)\|_{L^2(\O)\times L^2(\O)}=\|u\|_2+\|v\|_2.$$

\begin{theorem}\label{eh}
Let $M=\sup_\Gamma{T_w}$. For $\max \left\{C\|\nabla \Theta \|_\infty,M\right \}<\frac{\lambda}{C\|K\|_\infty}$ and small values of $\|K.\nabla \Theta\|_2$, there exists a unique solution $(\Psi,H)$ of the problem $(\ref{new1})$-$(\ref{bh2})$ in the space $\tilde W_\Psi \times \tilde W_H$.
\end{theorem}
\begin{proof}
Let us define the operator
\begin{equation*}
A\ :\  \tilde W_\Psi \times \tilde W_H \to L^2(\O) \times L^2(\O) 
\end{equation*}
such that $A(\Psi,H)=(A_1(\Psi,H),A_2(\Psi,H))$ with
\begin{align*}
A_1(\Psi,H)&=\Delta \Psi -K.\nabla H,\\ 
A_2(\Psi,H)&=\lambda \Delta H-\nabla H.(\nabla \Psi)^{\bot}-\nabla \Theta.(\nabla \Psi)^\bot.
\end{align*}
Let us remark that, using the Sobolev embedding theorem, we have  
$H^1(\Omega) \hookrightarrow L^4(\Omega)$ in such a way that $\nabla H.(\nabla \Psi)^{\bot}\in L^2(\O)$.
  
\noindent In term of the operator $A$, the equations $(\ref{new1})$-$(\ref{new2})$ can be rewritten as 
$$A(\Psi,H)=(K.\nabla \Theta,0).$$
Notice that $(\Psi,H)=(0,0)$ is a solution of $A(\Psi,H)=(0,0)$ and by the same argument as in Theorem \ref{unicite}, it is the only one. 

\noindent Now we want to show that the solution of $A(\Psi,H)=(K.\nabla \Theta,0)$ also exists for small values of $\|K.\nabla \Theta\|_2$. To this end, let us compute the Fr\'echet derivative of $A$. For $\phi \in \tilde W_\Psi$ and $G \in \tilde W_H$, we have
\begin{align*}
A(\phi,G)-(\Delta \phi-K.\nabla G,\lambda \Delta G-\nabla \Theta .(\nabla \phi)^\bot)
&:=A(\phi,G)-L(\phi,G)\\
&=(0,-\nabla G.(\nabla \phi)^\bot)\\
&=o(\|(\phi,G)\|_{\tilde W_\Psi \times \tilde W_H})
\end{align*}
because
\begin{align*}
\|(0,\nabla G.(\nabla \phi)^\bot)\|_{L^2(\O)\times L^2(\O)}&=\|\nabla G.(\nabla \phi)^\bot \|_{L^2(\O)}\\
&\leq \|\nabla G\|_{L^4(\O)}\|\nabla \phi\|_{L^4(\O)}\\
&\leq C_s^2 \|\nabla G\|_{H^1(\O)}\|\nabla \phi\|_{H^1(\O)}\\
&\leq 
C_s^2\|(\phi,G)\|^2_{\tilde W_\Psi \times \tilde W_H}
\end{align*}
where $C_s$ is the Sobolev constant corresponding to the continuity of the embedding $H^1(\Omega) \hookrightarrow L^4(\Omega)$.
Thus, $L$ defined by $L(\phi,G)=(\Delta \phi-K.\nabla G,\lambda \Delta G-\nabla \Theta.(\nabla \phi)^\bot)$ is the Fr\'echet derivative of $A$ at the point $(0,0)$, i.e.
$$A'(0,0).(\phi,G)=(\Delta \phi-K.\nabla G,\lambda \Delta G-\nabla \Theta. (\nabla \phi)^\bot).$$
For $f$ and $g$ in $L^2(\O)$ let us now consider the system $A'(0,0).(\phi,G)=(f,g)$ that can be written as
\begin{align}
-\Delta \phi+K.\nabla G&=f, \label{e1}\\
-\lambda \Delta G+\nabla \Theta. (\nabla \phi)^\bot&=g. \label{e2}
\end{align}
To prove the existence of a solution $(\Psi,H)$ of $(\ref{new1})$-$(\ref{bh2})$ it remains to show that the linear operator $A'(0,0):\tilde W_\Psi \times \tilde W_H \to L^2(\O) \times L^2(\O)$ is invertible. To this end, we must first prove that for every given $f$ and $g$ in $L^2(\O)$ the system (\ref{e1})-(\ref{e2}) admits at least a solution and secondly that for $(f,g)=(0,0)$ only $(\phi,G)=(0,0)$ is a solution of (\ref{e1})-(\ref{e2}).
\begin{itemize}

\item
First, we want to prove that for every given $f$ and $g$ in $L^2(\O)$ the system (\ref{e1})-(\ref{e2}) admits at least a solution. To this aim, let us define the operator $T=Q\circ S:G \mapsto G_1$ from $H^1(\O)$ into $H^1(\O)$ with $S: G\mapsto \phi$ where $\phi$ is the solution of
\begin{equation*}
-\Delta \phi+K.\nabla G=f \label{s1}
\end{equation*}
in $\Omega$ with the boundary conditions $\phi=0$ on $\Gamma_{1}$ and $\frac{\partial \phi}{\partial n}=0$ on $\Gamma_{2}$, and $Q: \phi \mapsto G_1$ where $G_1$ is the solution of
\begin{equation*}
-\lambda \Delta G_1+\nabla \Theta. (\nabla \phi)^\bot=g \label{q1}
\end{equation*}
in $\Omega$ with the boundary conditions $G_1=0$ on $\Gamma$.

Suppose now that $G$ and $G'$ are given in $H^1(\O)$. Let us consider 
$\phi=S(G)$, $\phi'=S(G')$ and 
$G_1=Q(\phi)$, $G_1'=Q(\phi')$. Setting $\bar G=G-G'$, $\bar \phi=\phi-\phi'$ and $\bar G_1=G_1-G_1'$, by (\ref{e1})-(\ref{e2}) we have the inequalities
\begin{align*}
\int_\O \|\nabla \bar \phi\|^2dx &=- \int_\O (K.\nabla \bar G)\phi dx
\leq C\|K\|_\infty \left ( \int_\O \|\nabla \bar \phi\| ^2dx\right )^\frac{1}{2}\left ( \int_\O \|\nabla \bar G\| ^2dx\right )^\frac{1}{2}
\end{align*}
and
\begin{align*}
\lambda \int_\O \|\nabla \bar G_1\|^2dx =-\int_\O \nabla \Theta.(\nabla \bar \phi)^\bot \bar G_1 dx \leq C\|\nabla \Theta \|_\infty \left ( \int_\O \|\nabla \bar \phi\| ^2dx\right )^\frac{1}{2}\left ( \int_\O \|\nabla \bar G_1\| ^2dx\right )^\frac{1}{2}.
\end{align*}
Combining these two inequalities, we obtain
$$\|\nabla \bar G_1\|_{L^2(\O)}\leq \frac{C^2\|K\|_\infty\|\nabla \Theta\|_\infty}{\lambda}\|\nabla \bar G\|_{L^2(\O)}$$
that shows us that if $$\frac{C^2\|K\|_\infty\|\nabla \Theta\|_\infty}{\lambda}<1$$
then $T$ is a contraction from $H^1(\O)$ into itself and admits a fixed point $G\in H^2(\O)$ that gives us a solution $(\phi,G)\in \tilde W_\Psi \times \tilde W_H$ of (\ref{e1})-(\ref{e2}). 

\item The system (\ref{e1})-(\ref{e2}) with $(f,g)=(0,0)$ admits $(0,0)$ for solution, let us show that this solution is unique. Let us suppose that $(\phi,G)\in \tilde W_\Psi \times \tilde W_H$ is a solution of (\ref{e1})-(\ref{e2}), multiplying (\ref{e1}) by $\phi$, (\ref{e2}) by $G$ and integrating on $\O$ leads to
$$\|\nabla G\|_2\leq \frac{C^2\|K\|_\infty\|\nabla \Theta\|_\infty}{\lambda}\|\nabla G\|_2$$
from which we deduce $G=0$ and $\phi=0$ if $C^2\|K\|_\infty\|\nabla \Theta\|_\infty<\lambda$.
\end{itemize}

This shows that, for small values of $\|K.\nabla \Theta\|_2$, the problem $A(\Psi,H)=(K.\nabla \Theta,0)$ does have solutions. Thus, for such values of $\Theta$ and $K$ and  $C^2\|K\|_\infty\|\nabla \Theta\|_\infty<\lambda$, the problem $(\ref{new1})$-$(\ref{bh2})$ admits at least one solution $(\phi,G)$ in $\tilde W_\Psi \times \tilde W_H$ and, as $\tilde W_\Psi \times \tilde W_H\subset W_\Psi \times W_H$, by Theorem \ref{unicite} it is unique if, in addition, we have $MC\|K\|_\infty<\lambda$.
\end{proof}

\begin{remark}
Since, in the previous Theorem we have
$$\|K.\nabla \Theta\|_2\leq \|K\|_\infty \|\nabla \Theta\|_\infty (\text{\rm mes }\O)^\frac{1}{2}$$
and 
$$\|\nabla \Theta \|_\infty<\frac{\lambda}{C^2\|K\|_\infty}$$
the condition $\|K.\nabla \Theta\|_2$ small is realized when $\frac{\lambda}{C^2}(\text{\rm mes }\O)^\frac{1}{2}$ is small. It is the case, for example, when the domain $\O$ is large and the parameter $\lambda$, that is the thermal diffusivity of the  porous medium, is small.
\end{remark}

\begin{corollary}
Let $T_w\in H^\frac{5}{2}(\Gamma)$ and $M=\sup_\Gamma{T_w}$. If $\max\left\{C\|\nabla \Theta \|_\infty,M\right \}<\frac{\lambda}{C\|K\|_\infty}$ there exists a unique solution $(\Psi,T)$ of the problem $(\ref{eq1})$-$(\ref{b2})$ in the space $\tilde W_\Psi \times H^2(\O)$ for small values of $\|K.\nabla \Theta\|_2$.
\end{corollary}
\begin{proof} It follows immediately from Theorem \ref{eh} and the fact that problems $(\ref{eq1})$-$(\ref{b2})$ and $(\ref{new1})$-$(\ref{bh2})$ are equivalent.
\end{proof}

\section{Conclusion}

In this paper, starting from the model of free convection in a fluid saturated porous medium near a semi-infinite vertical flat plate we have written an extension describing this phenomenon in a two-dimensional bounded domain. This new problem is given by two strongly coupled partial differential equations, that allows us to compute the stream function and the temperature of the fluid in the porous medium.

In a first approach of this complex problem, we have proved existence and uniqueness of a solution for small data when a part of the boundary of the domain is assumed to be impermeable.

\section*{Acknowledgement} 
The authors would like to thank Professor Herbert Amann for suggesting 
the idea for the existence proof and Professor Michel Chipot for his comments and helpful suggestions for proving Theorem 4.

\end{document}